\documentclass[reqno,12pt]{amsart}
\usepackage{amscd}
\usepackage{amssymb}
\usepackage[a4paper]{geometry}
\usepackage{amstext}
\usepackage{mathrsfs}
\usepackage{times} 
\usepackage{todonotes}
\usepackage[final]{hyperref}
\hypersetup{unicode= false, colorlinks=true, linkcolor=blue,
anchorcolor=blue, citecolor=red, filecolor=red, menucolor=blue,
pagecolor=blue, urlcolor=blue}
\numberwithin{equation}{section}
\usepackage[none]{hyphenat}

\def\a{\alpha}

\def\bege{\begin{equation}} \def\ende{\end{equation}}
  \def\a{\alpha} 
  \def\ve{\varepsilon} 
 
\def\begr{\begin{eqnarray}} \def\endr{\end{eqnarray}}
\def\qand{\quad\mbox{ and }\quad}

\def\bege{\begin{equation}} \def\ende{\end{equation}}
\def\begr{\begin{eqnarray}} \def\endr{\end{eqnarray}}
\def\bnum{\begin{enumerate}} \def\enum{\end{enumerate}}

\newcommand{\NN}{\mathbb{N}}

\newcommand{\CC}{\mathbb{C}}

\DeclareGraphicsExtensions{.pdf,.png,.jpg}
\newtheorem{theorem}{Theorem}[section]

\newtheorem{lemma}{Lemma}[section]
\newtheorem{corollary}{Corollary}[section]

\title[Generalized Volterra companion operators on Fock spaces] {Generalized Volterra companion operators on Fock spaces }
\author [Tesfa  Mengestie]{Tesfa  Mengestie }
\address{Department of Mathematical Sciences \\
Stord/Haugesund University College (HSH)\\
Klingenbergvegen 8, N-5414 Stord, Norway}
\email{tesfantnu@gmail.com}


\begin{document}

\begin{abstract}
We characterize the bounded and compact generalized Volterra companion integral operators on Fock spaces
 acting between the standard Fock spaces. As a special case, we
prove that there exist no nontrivial compact Volterra companion
integral and multiplication operators  on Fock spaces.  We also obtain  asymptotic
estimates for the norm of these operators. 
\end{abstract}

\maketitle

\section{Introduction}

For  holomorphic functions $f$ and $g$, the  Volterra type
integral operator $V_g$ and its companion $J_g$  are defined by
\begin{align*}
 V_gf(z)= \int_0^z f(w)g'(w) dw \ \ \ \text{and} \ \  \ J_gf(z)= \int_0^z f'(w)g(w) dw.
\end{align*}
Applying integration by parts in any one of the above integrals
gives the relation $$V_g f+ J_g f= M_g f-f(0)g(0),$$ where $M_g
f= g f$ is the multiplication operator of symbol $g$.   These
integral type operators have been studied extensively on  various spaces of
analytic functions with the aim to explore the connection between
their operator theoretic behaviours with the function theoretic properties of  the
symbols $g$ especially after the works of Pommerenke \cite{Pom}, and
 Aleman and Siskakis \cite{Alsi1,Alsi2} on
Hardy and Bergman spaces.   For more information on the subject,
we refer to \cite{Alman,ALC,Si} and the related references
therein.

The idea  to extend the operators $V_g$ and $J_g$ was first
raised by S. Li and S. Stevi\'c in 2006. They  eventually proposed introducing the following
 operators induced by pairs of holomorphic symbols
$(g,\psi)$:
\begin{align}
\label{OKK}
V_g^{\psi} f(z)= \int_0^z f(\psi(w)) g'(w) dw, \ \ \  \ \
\ C_g^\psi f(z)= \int_0^{\psi(z)} f(w) g'(w)dw, \\
J_{(g,\psi)} f(z)= \int_0^z f'(\psi(w)) g(w) dw, \ \text{and} \  \ C_{(g,\psi)} f(z)= \int_0^{\psi(z)} f'(w) g(w)dw,
\label{OKJ}
\end{align} and studied their operator theoretic properties in terms of the pairs  $(g,\psi)$  on some spaces of analytic
functions on the unit disk \cite{ LIS, jmaa345}. For more  recent
results on the class of  operators in  \eqref{OKK}, one may consult the materials for instance  in  \cite{TM, TM0, UK2}.

Studying
operators $V_g^{\psi}$ and $C_g^\psi$ attracted somewhat more
attention partly due to  the fact that their bounded and
compact properties are related to the notion of Carleson measures,
which are readily available for several known spaces. In contrast,
relatively little is known on the operators $J_{(g,\psi)}$ and
$C_{(g,\psi)}$ apart from some cases where the target space is
restricted to be a growth type space with norm defined  in terms
of derivatives. We may mention that these class of operators  have
also found applications in the study of linear isometries of
spaces of analytic functions.  A prototype example in this case is the following. Let  $S^p$ denotes the
space of all analytic functions $f$ in the unit disc for which its
derivative $f'$ belongs to the Hardy space $H^p.$  Then it has been show that   for
$p\neq2,$ any surjective isometry $T$ of $S^p$ under the norm
$\|f\|_{S^p}= |f(0)|+ \|f'\|_{H^p}$ is of the form $$Tf= \lambda
f(0)+ \lambda J_{(g,\psi)}f$$ for some unimodular $\lambda$ in $ \CC$,
 a nonconstant inner function $\psi$ and a function $g$ in $H^p$ \cite{FJ}.

The bounded and compact properties of the  class of operators in \eqref{OKK}  when acting between   Fock spaces were studied  in
\cite{TM,TM0}. In this note, we  continue that  line of research
for the remaining class of   operators $J_{(g,\psi)}
$ and $C_{(\psi,g)}$ which are called the generalized Volterra companion
integral operators. Observe that the particular  choice $\psi(z)= z$
reduces  both $J_{(g,\psi)} $ and $C_{(\psi,g)}$  to the operator
$J_g$. On the other hand, setting
$ g=\psi'$  and $g = 1$ respectively  reduce  the operators $J_{(g,\psi)} $ and $C_{(\psi,g)}$  to the
composition operator $ C_\psi$  up to a constant.   As will be seen latter, a consequence of our main results
shows that there exist no nontrivial compact Volterra
companion integral and multiplication operators acting between Fock spaces.

The classical weighted Fock space $\mathcal{F}_\alpha^p$ consists
of all entire functions $f$ for which
\begin{equation*}
\|f\|_{p}^p=  \frac{\alpha p}{2\pi}\int_{\CC} |f(z)|^p
e^{-\frac{\alpha p}{2}|z|^2} dm(z) <\infty,
\end{equation*} where  $0 < p <\infty,$ $ \alpha >0 $ is a parameter, and  $dm$ denotes the
usual Lebesgue area  measure on $\CC$. For  $p= \infty,$ the
growth type space $\mathcal{F}_\alpha^\infty$  contains all entire
functions $f$ such that
\begin{equation*}
\|f\|_{\infty}= \sup_{z\in \CC}
|f(z)|e^{-\frac{\alpha}{2}|z|^2} <\infty.
\end{equation*} The space  $\mathcal{F}_\alpha^2$, in particular, is  a reproducing kernel Hilbert space with
kernel function $K_{w}(z)= e^{\alpha \langle z, w\rangle}$ and
normalized kernel function $k_{w}(z)= e^{
\alpha\langle z, w\rangle}-\alpha|w|^2/2.$

We close this introduction with a word on notation: The notation $U(z)\lesssim V(z)$ (or
equivalently $V(z)\gtrsim U(z)$) means that there is a constant
$C$ such that $U(z)\leq CV(z)$ holds for all $z$ in the set of a
question. We write $U(z)\simeq V(z)$ if both $U(z)\lesssim V(z)$
and $V(z)\lesssim U(z)$. In addition, we denote by $L^p$ the Lebesque spaces
$L^p(\CC,dm)$.

\section{The main results}
In this section, we will present our main results. We may first set  $$P_\psi (z)=
\frac{e^{\frac{\alpha}{2}(|\psi(z)|^2-|z|^2)}}{1+|z|}\qand
Q_g(z)=\frac{|g(z)|e^{-\frac{\alpha }{2}|z|^2}}{1+|z|}.$$
 Then our
results are expressed in terms of the functions $$
M_{(g,\psi)}(z)= |g(z)|\big(|\psi(z)|+1\big) P_\psi (z) \qand
M_{(g(\psi),\psi)}(z)= |g(\psi(z))|\big(|\psi(z)|+1\big) P_\psi (z), $$ and
generalized Berezin type integral transforms:
\begin{align*}
  B_{(|g(\psi)|^q,\psi)} (w)=\int_{\CC} \big|k_w(\psi(z)\big)\big|^q \Big(\big(|w|+1)|\psi'(z)| Q_{g(\psi)}(z)\Big)^q dm(z)
  \end{align*}
  \and
  \begin{align*}
B_{(|g|^q,\psi)} (w)= \int_{\CC} \big|k_w(\psi(z))\big|^q
\Big( \big(|w|+1\big)Q_g(z)\Big)^q dm(z). \ \ \ \ \ \ \ \ \ \ \ \
\end{align*}
We now state our first main result.
\begin{theorem}\label{thm1}
Let $0<p\leq q \leq\infty$ and $(g,\psi)$ be a pair of entire
functions. Then
\begin{enumerate}
 \item $J_{(g,\psi)}: \mathcal{F}_\alpha^p \to \mathcal{F}_\alpha^q$ is bounded if and only if
 $M_{(g,\psi)}$ when $q= \infty$ and  $B_{(|g|^q,\psi)} $ when $q<\infty$  belong to $L^\infty$. In this case,
 we also have
     \begin{equation}
     \label{bdd1}
    \|J_{(g,\psi)}\|\simeq  \begin{cases}
    \|M_{(g,\psi)}\|_{L^\infty}, & q= \infty\\
    \|B_{(|g|^q,\psi)}\|_{L^\infty}^{1/q},   &q<\infty.
    \end{cases}
     \end{equation}
 \item  $J_{(g,\psi)}: \mathcal{F}_\alpha^p \to \mathcal{F}_\alpha^q$ is compact if and only if
 \begin{center}
 \begin{enumerate}
 \item $\lim_{|z|\to \infty}B_{(|g|^q,\psi)}(z)=0$  for
$q<\infty$;
\item   it is bounded and $ \lim_{|\psi(z)|\to \infty
}M_{(g,\psi)}(z)= 0$ when  $ q=\infty$ .
\end{enumerate}
\end{center}

\item $C_{(g,\psi)}: \mathcal{F}_\alpha^p \to
\mathcal{F}_\alpha^q$ is  bounded if and only if
     $M_{(g(\psi),\psi)}$ when $q= \infty$ and $B_{(|g(\psi)|^q,\psi)} $ when $q<\infty$ belong to $ L^\infty$.
     We also estimate the norm by
     \begin{equation}
     \label{bdd2}
    \|C_{(g,\psi)}\|\simeq  \begin{cases}
    \|M_{(g(\psi),\psi)}\|_{L^\infty}, & q= \infty\\
    \|B_{(|g(\psi)|^q,\psi)}\|_{L^\infty}^{1/q},   &q<\infty.
    \end{cases}
\end{equation}
\item $C_{(g,\psi)}: \mathcal{F}_\alpha^p \to
\mathcal{F}_\alpha^q$ is compact if and only if
\begin{center}
\begin{enumerate}
\item
$\lim_{|z|\to\infty}B_{(|g(\psi)|^q,\psi)}(z)=0$ for $ q<\infty$;
\item it is bounded and  $ \lim_{|\psi(z)|\to \infty
}M_{(g(\psi),\psi)}(z)= 0$ when $q=\infty$.
\end{enumerate}
\end{center}
\end{enumerate}
\end{theorem}
It may be noted that the conditions on the  preceding results  do not depend on  exponent $p$ from
 the  domain space $\mathcal{F}_\alpha^p$  apart from the sole  assumption that  $p\leq q,$ where $q$ is the exponent on the target space. It means that if there exists a $p_o \leq q $ for which the map $ J_{(g,\psi)}$ or $C_{(g,\psi)}$ is  bounded (compact) from $\mathcal{F}_\alpha^{p_o}$ to $\mathcal{F}_\alpha^{q}$, the same conclusion holds when we replace the domain space by $\mathcal{F}_\alpha^{p}$ for any $p\leq q.$ A similar phenomena  was observed in \cite{TM,TM0} with  the class of operators in \eqref{OKK} .

Recall  that since the Fock spaces are nested, $\mathcal{F}_\alpha^p\subseteq \mathcal{F}_\alpha^p$ whenever
$p\leq q$ \cite{SJR},  the conditions in Theorem~\ref{thm1} are from mapping smaller spaces into larger spaces under the maps $J_{(g,\psi)}$ and $C_{(g,\psi)}$. Conversely, when  we map larger spaces into smaller spaces with the same mappings, we get the following stronger  integrability conditions.

 \begin{theorem}\label{thm2}
 Let $0<q<p \leq\infty$ and $(g,\psi)$ be a pair of entire functions. Then
 \begin{enumerate}
 \item $J_{(g,\psi)}: \mathcal{F}_\alpha^p \to \mathcal{F}_\alpha^q$ is bounded (compact) if and only if
 $B_{(|g|^q,\psi)}$ belongs to $L^{1}$ when  $p=\infty$ and  to $L^{p/(p-q)}$ whenever
 $p<\infty.$  Furthermore, we have
      \begin{equation}
     \label{bdd3}
    \|J_{(g,\psi)}\|\simeq  \begin{cases}
   \|B_{(|g|^q,\psi)}\|_{L^1},  &\ p=\infty\\
   \|B_{(|g|^q,\psi)}\|_{L^{\frac{p}{p-q}}},&  p<\infty.
    \end{cases}
     \end{equation}
  \item $C_{(g,\psi)}: \mathcal{F}_\alpha^p \to \mathcal{F}_\alpha^q$ is  bounded (compact) if and only if
  $B_{(|g\psi|^q,\psi)}$ belongs to $L^{1}$ for $p=\infty$ and  to $L^{p/(p-q)}$ whenever
 $p<\infty.$ Furthermore, we have
     \begin{equation}
     \label{bdd4}
    \|J_{(g,\psi)}\|\simeq  \begin{cases}
   \|B_{(|g(\psi)|^q,\psi)}\|_{L^1},  &\ p=\infty\\
   \|B_{(|g(\psi)|^q,\psi)}\|_{L^{\frac{p}{p-q}}},&  p<\infty.
    \end{cases}
     \end{equation}
  \end{enumerate}
\end{theorem}
Unlike the case for Theorem~\ref{thm1}, the conditions in the preceding theorem
rely on both the domain and the targets space exponents $p$ and $q$ except for  the case when
$p= \infty$. In this case, the corresponding condition is independent of the target space exponent $q$ as long as
$q<p.$

From the relation $V_g f+ J_g f= M_g f-f(0) g(0),$  we observe
that if any two of the operators $V_g, J_g$ and $M_g$ are bounded
so does the third onve. Interestingly, in  Fock spaces  more can be
said, namely that $M_g$ is bounded (compact) if and only if  so is
the operator $J_g$. We formulate this observation as a corollary below and prove our assertion in the next section.

\begin{corollary}\label{cor1}
Let $g$ be an entire function on $\CC$. Then
\begin{enumerate}
\item if    $0<p\leq q\leq \infty,$ then  the operator $J_g: \mathcal{F}_\alpha^p \to \mathcal{F}_\alpha^q$ is bounded (respect. compact) if and only if so is $M_g: \mathcal{F}_\alpha^p \to \mathcal{F}_\alpha^q$,
 and this holds if and only if  g is a constant (respect. zero) function. 
  \item if $ 0<q<p\leq  \infty,  $ then the operator $J_g: \mathcal{F}_\alpha^p \to \mathcal{F}_\alpha^q$ is   bounded or compact if and only if so is $M_g: \mathcal{F}_\alpha^p \to \mathcal{F}_\alpha^q$, and this holds if and only if  $g$ is the zero function.
\end{enumerate}
 \end{corollary}

The results in the corollary verify that there exist no nonzero
compact Volterra companion integral  and multiplication operators acting
between Fock spaces. Furthermore, it  has now become clear that
more symbols $g$ are admissible in inducing  bounded or compact
$V_g$ than $J_g$. More specifically, a number of results from
\cite{Olivia,TM,TM0} ensure  that $V_g:\mathcal{F}_\alpha^p \to
\mathcal{F}_\alpha^q, 0 <p\leq q\leq\infty,$ is bounded if and
only if $g$ is a complex polynomial of degree not exceeding $2,$
and its compactness  holds if and only if its   degree does not
exceed $1.$  On the other hand, if $p>q,$ then $V_g$ is bounded or compact if
and only if $g$ is again  a  polynomial of  degree not exceeding $1.$ It means that the  bounded (compact) properties of
the sum  $V_g+ J_g$  merely depends on the  boundedness
(compactness) of  the summand  $J_g$ when all the three operators
act between  Fock spaces.

By simply scaling  $z$ as $\beta z$ for some $ |\beta|<1$,  it is easily seen that a wider class of  symbols $g$ are  admissible in giving
rise to bounded (compact) $J_{(g,\psi)}$ than  those $g$ guaranteed in the corollary. For instance  if we  set $\psi(z)= \beta |s|, |\beta|<1$ and $g$ be entire function such that $|g(z)|^2 \leq e^{\gamma|z|^2}$ where $ |\beta|^2 + \gamma <1.$ Then the pair $(g,\psi)$ satisfies the condition in part (i) of Theorem~\ref{thm1}.
This gives another impetus for the need to  take further the study of Volterra companion integral  and composition operators to the generalized cases $J_{(g,\psi)}$ and $C_{(\psi,g)}$.\\

Combining the results from Theorem~\ref{thm1}, Theorem~\ref{thm2} and the  corresponding boundedness (compactness) results from \cite{TM,TM0}, we immediately deduce the following.
 \begin{corollary}\label{cor2}
 Let $g$ be a holomorphic function on $\CC$ and $0<p,q \leq \infty.$ Then if
 \begin{enumerate}
 \item  $J_{(g,\psi)}: \mathcal{F}_\alpha^p \to \mathcal{F}_\alpha^q$ is bounded (compact) so is the map
 $V_g^\psi: \mathcal{F}_\alpha^p \to \mathcal{F}_\alpha^q$.
 \item  $C_{(g,\psi)}: \mathcal{F}_\alpha^p \to \mathcal{F}_\alpha^q$ is bounded (compact) so is the map
 $C_g^\psi: \mathcal{F}_\alpha^p \to \mathcal{F}_\alpha^q$.
 \end{enumerate}
  \end{corollary}
  The converses of the statements in the corollary in general fail. To find a  simple counterexample, we may simply
  set $\psi(z)= z.$ Then the class of operators in \eqref{OKK} reduces to the Volterra type integral operator
  $V_g$ while those at \eqref{OKJ} reduce to its companion operator $J_g.$ Then the desired conclusion follows from the analysis in the  paragraph immediately after Corollary~\ref{cor1} above.

\section{Auxiliary Results}
In this section we collect some  auxiliary results that will be used  in our
subsequent considerations. Our first lemma provides a  criteria for
compactness of $J_{(g,\psi)}$ and $C_{(g,\psi)}$  when acting between  Fock spaces.

\begin{lemma}\label{complem}
Let $0<q,p \leq\infty$ and $(g,\psi)$ be a pair of entire
functions. Then
\begin{enumerate}
\item $J_{(g,\psi)}: \mathcal{F}_\alpha^p \to
\mathcal{F}_\alpha^q $ is compact if and only if $\|
J_{(g,\psi)}f_n\|_{q} \to 0$ as $n\to \infty$ for each uniformly bounded
sequence $(f_n)_{n\in \NN}$ in  $\mathcal{F}_\alpha^p$ converging
to zero uniformly on compact subsets of $\CC$ as $n\to \infty.$
\item  A similar statement holds when we replace the operator  $J_{(g,\psi)}$ by $C_{(g,\psi)}$ in
(i).
\end{enumerate}
\end{lemma}
The lemma can be proved following standard arguments, and will  be used repeatedly in what follows  without mentioning
it over and over again.

For  $q\geq 0$, we set $\varphi_q(z)= (1+|z|)^q$ to  be a weight function on $\CC$. Then for each $p\geq 1,$ we introduce weighted $L_{\varphi_q}^p$ spaces
  consisting of  all measurable functions $f$ on $\CC$ such that
  \begin{equation*}
  \|f\|_{L_{\varphi_q}^p}^p= \int_{\CC} |f(z)\varphi_q(z)|^p dm(z) <\infty
  \end{equation*} for  finite $ p$  and  when $p=\infty,$ the corresponding norm is given by
  \begin{equation*}
  \|f\|_{L_{\varphi_q}^\infty}= \sup_{z\in \CC}  |f(z)\varphi_q(z)|<\infty.
  \end{equation*}
 For a Borel measure $\mu$ on $\CC$, we also define a   Berezin type integral
  transform associated with  it by
  \begin{equation}
  \widetilde{\mu}_q (z)= \int_{\CC} (1+|z|)^q |k_{z}(\zeta)|^q e^{-\frac{\alpha q}{2}|\zeta|^2} d\mu(\zeta).
  \end{equation}
In particular  when  $\mu$  is a measure such that  $d\mu (z)=f(z)dm(z)$ for a given
measurable function $f$, we prove the following.

  \begin{lemma}\label{gen} let $1\leq p\leq \infty$, $r>0$ and $0\leq q<\infty$. Then  the operators
  $f \mapsto f_r$ and  $f \mapsto \widetilde{f}_q$  from  $L_{\varphi_q}^p$ to $L^p$ are bounded where $f_r(z)= (1+|z|)^q \mu(D(z,r)),$  $D(z,r)= \{\zeta \in \CC: |z-\zeta|<r\}$, and $$\widetilde{f}_q(z)= \int_{\CC} (1+|z|)^q |k_{z}(\zeta)|^q e^{-\frac{\alpha q}{2}|\zeta|^2} d\mu(\zeta).$$
    \end{lemma}
 \begin{proof} We mention that for the case when $q=0$, the lemma was first proved in \cite{ZHXL}. We use  interpolation arguments  between the  Lebesgue spaces  $L_{w}^p$ and  $L^p$  and extend the techniques there. Thus, it suffices to establish the statements for $p=1$ and $p=\infty.$ We begin with the case $p=1$ and  apply Fubini's theorem to estimate
  \begin{align*}
  \|\widetilde{f}_q\|_{L^1}&= \int_{\CC} \bigg|\int_{\CC}  |(1+|z|)^q |k_z(w)|^q e^{-\frac{\alpha q}{2} |w|^2}d\mu(w) \bigg| dm(z)\nonumber\\
  &\leq \int_{\CC} \bigg(\int_{\CC} |k_z(w)|^q e^{-\frac{\alpha q}{2} |w|^2} (1+|z|)^q dm(z) \bigg)  |f(w)|dm(w)\nonumber\\
   &=\int_{\CC} \bigg(\int_{\CC} (1+|z|)^q e^{-\frac{\alpha q}{2} |z-w|^2} dm(z) \bigg)  |f(w)|dm(w)\nonumber\\
   &\lesssim \int_{\CC}  (1+|w|)^q |f(w)|dm(w) = \|f\|_{L_{\varphi_q}^1}.
  \end{align*} On the other hand, applying Fubini's theorem  again and the fact that
  $\chi_{D(\zeta,r)}(z)= \chi_{D(z,r)}(\zeta)$ for all $\zeta$ and $z$ in $\CC$, we have
  \begin{align*}
  \|f_r\|_{L^1}=\int_{\CC} (1+|z|)^q \mu(D(z,r))dm(z) \leq  \int_{\CC} (1+|z|)^q \int_{D(z,r)} |f(\zeta)|dm(\zeta) dm(z)\nonumber\\
  = \int_{\CC} |f(\zeta)| \int_{\CC} \chi_{D(\zeta,r)}(z) (1+|z|)^q dm(z) dm(\zeta)\nonumber\\
  \simeq \int_{\CC} |f(\zeta)| (1+|\zeta|)^q dm(\zeta)= \ \|f\|_{L_{\varphi_q}^1},
  \end{align*} where the last estimate  follows  since  $1+|z| \simeq 1+|\zeta|$ for each
  $z$ in $D(\zeta,r)$.\\
  We now proceed to show the case for  $p= \infty.$  For each $f\in L^\infty_{\varphi_q}$,  it easily follows that
  \begin{align*}
  \sup_{z\in \CC}|f_r(z)|\leq \sup_{z\in \CC} (1+|z|)^q \int_{D(z,r)} |f(\zeta)|dm(\zeta)\simeq \sup_{z\in \CC}\int_{D(z,r)} (1+|\zeta|)^q |f(\zeta)|dm(\zeta)\\
  \lesssim \sup_{z\in \CC} \sup_{\zeta \in D(z,r)} (1+|\zeta|)^q |f(\zeta)|\leq \|f\|_{L_{\varphi_q}^\infty}.
  \end{align*}
  Seemingly,    we also have
  \begin{align}
  \label{del}
  \sup_{z\in \CC}|\widetilde{f}_q(z)|\leq \sup_{z\in \CC} (1+|z|)^q \int_{\CC}|k_z(\zeta)|^q e^{-\frac{q\alpha}{2}|\zeta|^2} |f(\zeta)|dm(\zeta)\nonumber\\
  =  \ \ \ \  \sup_{z\in \CC} (1+|z|)^q \int_{\CC} e^{-\frac{q\alpha}{2}|z-\zeta|^2} |f(\zeta)|dm(\zeta)\nonumber\\
   = \sup_{z\in \CC} (1+|z|)^q \int_{\CC} e^{-\frac{q\alpha}{2}|\tau|^2} |f(z-\tau)|dm(\tau).
  \end{align}
Now observe that since  $f$ belongs to $ L_{\varphi_q}^\infty$, it satisfies $  |f(\zeta)| \lesssim (1+ |\zeta|)^{-q}$
 for each $\zeta$ in $\CC$.  This means that the last integral in \eqref{del} is uniformly
 bounded independent of $z.$  If we set  $M_1$ to be one of such a bound and
 $$ M= \sup_{f\in  L_{\varphi_q}^\infty } \sup_{z\in \CC}|f(z)| <\infty,$$  then we see that the integral is indeed
bounded by  $M_1 M^2 |f(z)|$ for any $z$. Taking this into account, we  obtain
\begin{align*}
\sup_{z\in \CC}|\widetilde{f}_q(z)| \leq  \sup_{z\in \CC} (1+|z|)^q \int_{\CC} e^{-\frac{q\alpha}{2}|\tau|^2} |f(z-\tau)|dm(\tau)\\
\leq \sup_{z\in \CC} (1+|z|)^q M_1M^2 |f(z)|\simeq \|f\|_{L_{\varphi_q}^\infty}
\end{align*} and completes the proof.
\end{proof}
\begin{lemma}
\label{forall}
Let $\mu$ be a nonnegative measure on $\CC$, \ \   $D_{rq\mu}(z)=(1+|z|)^q \mu\big(D(z, r)\big)$
 for positive values $r$ and $q$,  and  $ 0< p\leq \infty$.
 Then if $D_{\delta q\mu}$ belongs to $L^p$ for some
$\delta >0,$ then $D_{rq\mu}$ belongs to $L^p$ for all $r>0.$
\end{lemma}
\begin{proof}
For each $\tau$ in $\CC$ we may write
\begin{align*}
\int_{D(\tau, r)} (1+|z|)^q \mu(D(z,\delta)) dm(z)= \int_{\CC} \int_{\CC}   (1+|z|)^q \chi_{D(\tau,r)}(z)  \chi_{D(z,\delta)}(\zeta) d\mu(\zeta) dm(z).
\end{align*} Using again the simple fact that  $\chi_{D(z,\delta)}(\zeta)= \chi_{D(\zeta,\delta)}(z)$, the double integral
above is easily seen to be equal to
\begin{align*}
\int_{\CC} \int_{D(\zeta,\delta) \cap D(\tau, r)} (1+|z|)^q dm(z) d\mu(\zeta)\simeq \int_{\CC} (1+|\tau|)^q m\big(D(\zeta,\delta) \cap D(\tau, r)\big) d\mu(\zeta)\\
\geq  (1+|\tau|)^q \int_{D(\tau, r)}m\big(D(\zeta,\delta) \cap D(\tau, r)\big) d\mu(\zeta),
\end{align*} where $m(E)$ refers to the Lebesque area measure of set $E.$   Clearly,
the right hand quantity is bounded from below by
\begin{equation*}
 (1+|\tau|)^q \mu(D(\tau, r)) \inf_{\zeta\in D(\tau, r)} m\big(D(\zeta,\delta) \cap D(\tau, r)\big) \gtrsim (1+|\tau|)^q \mu(D(\tau, r)),
\end{equation*} where the lower estimate follows  since $\zeta\in D(\tau, r)$, there obviously  exists a disc $D(\tau_0, r_0)$ contained in $D(\zeta,\delta) \cap D(\tau, r)$  with $m\big( D(\tau_0, r_0)\big)= \pi r_0^2.$ From the above analysis, we conclude
\begin{equation}
\label{ess3}
(1+|\tau|)^q \mu(D(\tau, r)) \lesssim \int_{D(\tau, r)} (1+|z|)^q \mu(D(z,\delta)) dm(z).
\end{equation} If we now  set $f(z)=\mu(D(z,\delta)),  $ then the estimate above  along with Lemma~\ref{gen} ensure
that \begin{equation}\label{forbig}
\|D_{rq\mu}\|_{L^p} \lesssim \|f_\delta\|_{L^p} \lesssim\|f\|_{L_{\varphi_q}^p}=\|D_{\delta q\mu}\|_{L^p} <\infty
\end{equation} for each $p\geq1$ and any $r>0.$ We need to deduce the same when  $0<p<1.$  In this case,  \eqref{ess3} and Fubini's theorem  imply
\begin{align*}
\int_{\CC}(1+|\tau|)^{pq} (\mu(D(\tau, r)))^p dm(\tau)&\lesssim \int_{\CC}\bigg(\int_{D(\tau, r)} (1+|z|)^q \mu(D(z,\delta)) dm(z)\bigg)^pdm(\tau)\\
&\leq \int_{\CC}\int_{D(\tau, r)}  (1+|z|)^{pq} (\mu(D(z,\delta)))^p dm(z)dm(\tau)\\
&=\int_{\CC} \int_{\CC} \chi_{D(\tau,r)}(z)  \frac{(\mu(D(z,\delta)))^p}{(1+|z|)^{-pq}} dm(z)dm(\tau)\\
&=\int_{\CC}\bigg( \int_{\CC} \chi_{D(z,r)}(\tau)dm(\tau)\bigg)  \frac{(\mu(D(z,\delta)))^p}{(1+|z|)^{-pq}} dm(z)\\
&= \pi r^2\int_{\CC} (1+|z|)^{pq}(\mu(D(z,\delta)))^p dm(z)<\infty
\end{align*} independent of the choice of $r$ again.  This along with  \eqref{forbig} establishes
our claim for all exponent $p.$
\end{proof}
 We  next   recall  the notion of
lattice for the complex plane $\CC$. For a positive $r$, we  say that a sequence of distinct points
$(z_k)_{k\in \NN} \subset \CC$ is  an $r/2-$ lattice  for  $\CC$ if the  sequence of the discs $D(z_k, r), \ k\in \NN $ constitutes a covering of $\CC$ and the discs $ D(z_k, r/2)$ are mutually disjoint. An interesting example of such a lattice  can be found in \cite{ZHXL}.
\begin{lemma} \label{covering} Let $r>0$ and $(z_k)_{k\in \NN}$ be  an $r/2-$ lattice  for  $\CC$. Then there exists a positive integer $N_{\max}$ such that every point in $\CC$ belongs to at
most $N_{\max}$ of the discs $D(z_k, 2r)$.
\end{lemma}
The proof of the lemma can be found in \cite{KH,KZH} where in
\cite{KH} a more general setting has been considered.  The
sequence $z_k, k\in \NN$ will refer to such a fixed  $r/2$ lattice in the
remaining part of the paper. We now establish a basic lemma that
will be used in the proof of the necessity parts  of  our main results.
The lemma is also of its own interest.

\begin{lemma} \label{basic}
Let $\mu$ be a nonnegative measure on $\CC, \ q\geq0$ and $ 0< p\leq
\infty$. Then the following statements are equivalent.
\begin{enumerate}
\item The function $\widetilde{\mu}_q$ belongs to $L^p$;
\item  $D_{r q\mu}$ belongs to $L^p$ for
some (or any) $r>0$;
\item The sequence
$\big((|z_j|+1)^q\mu(D(z_j,r))\big)_{j\in \NN}$ belongs to $\ell^p$ for
some (or any) $r>0.$ Moreover, we have
\begin{equation}
\label{normre1} \|\widetilde \mu_q\|_{L^p}\simeq \|D_{rq\mu}\|_{L^p} \simeq
\|\big((1+|z_j|)^q\mu(D(z_j,r))\big)\|_{l^p}.\end{equation}
\end{enumerate}
\end{lemma}
\begin{proof} We begin with the proof of (i) implies (ii).  For any
nonnegative $r$ and $q$, we have
\begin{align*}
\label{firstnec}
D_{rq \mu}(z) =&(|z|+1)^q\int_{D(z,r)}
d\mu(w)\nonumber\\
\lesssim& \int_{D(z,r)}(|z|+1)^q e^{\frac{\alpha
q}{2}(r^2-|w-z|^2)}d\mu(w)\lesssim \tilde{\mu}_q(z)
\end{align*}
from which we  obtain the estimate
\begin{equation}
\label{normest1} \|D_{rq\mu} \|_{L^p} \lesssim
\|\widetilde{\mu}_q\|_{L^p}<\infty
\end{equation} independently of  the choices of $r$ and exponent $p$.\\
  We now prove  the equivalency of  statements (ii) and (iii), and assume that  part (ii) holds. Note that by the triangle inequality
       \begin{align}
       \label{tri}
       \mu(D(z,2r)) \geq \mu(D(z_j,r))
       \end{align}for each $z$ in $D(z_j,r)$. This  and Lemma~\ref{forall} imply
       \begin{align*}
    \infty> N_{\max}\int_{\CC}\bigg(\frac{\mu(D(z,2r))}{(|z|+1)^{-q}}\bigg)^p dm(z) \geq \sum_{j=1}^\infty \int_{D(z_j,r)} \bigg(\frac{\mu(D(z,2r))}{(|z|+1)^{-q}}\bigg)^p dm(z)\nonumber\\
       \geq \sum_{j=1}^\infty \int_{D(z_j,r)} \big((|z|+1)^q\mu(D(z_j,r))\big)^pdm(z)
       \geq \frac{\pi r^2}{2}  \sum_{j=1}^\infty \bigg(\frac{\mu(D(z,2r))}{(|z|+1)^{-q}}\bigg)^p.
      \end{align*} Since for each $z$, there exists $j$ for which $z$ belongs to the disc $D(z_j,r)$, the case $p= \infty$
      follows easily from the relation in \eqref{tri}. Therefore,
      \begin{equation}
      \label{est1}
      \|\big((|z_j|+1)^q\mu(D(z_j,r))\big)\|_{\ell^p}\lesssim \|D_{2rq\mu}\|_{L^p}
      \lesssim \|D_{rq\mu}\|_{L^p}
            \end{equation} for all positive $r$ and all $p.$

 Conversely, assume that  the sequence $\big( \mu(D(z_j,r)) (1+|z_j|)^q\big)_{j\in \NN}$ belongs to $\ell^p$. Then
            \begin{align}
            \label{normest2}
            \int_{\CC}D_{rq\mu}^p (z)dm(z)\simeq \sum_{j=1} \int_{D(z_j,r)} |\mu(D(z,r))|^p (1+|z|)^{qp} dm(z)\nonumber\\
            \leq \sum_{j=1} \int_{D(z_j,r)} \big(\mu(D(z_j,r))\big)^p (1+|z|)^{qp} dm(z)\nonumber\\
            \simeq \sum_{j=1} \int_{D(z_j,r)} \big((\mu(D(z_j,r))\big)^p (1+|z_j|)^{qp}.
            \end{align}The case $p=\infty$  holds  trivially again. Thus, from \eqref{normest2} we establish the estimate
            \begin{equation}
            \label{ess2}
            \|D_{rq\mu} \|_{L^p} \lesssim \|\big((|z_j|+1)^q\mu(D(z_j,r))\big)\|_{\ell^p}.
            \end{equation}
To complete the proof of the lemma, we next show that (ii) implies (i).
By Lemma~2.1 of \cite{ZHXL} applied to the function $k_w$,  we have
\begin{align}
|k_w(\zeta)|^q  e^{-\frac{q\alpha}{2}|\zeta|^2} \lesssim \int_{D(\zeta,r)}  |k_w(z)|^q e^{-\frac{q\alpha}{2}|z|^2}dm(z).
\end{align} Multiplying both sides of the inequality by $\varphi_q(w)$ and subsequently integrating against the measure $\mu$ give
\begin{align}
\widetilde{\mu}_q (w)&= \int_{\CC} (|w|+1)^q|k_w(\zeta)|^qe^{-\frac{q\alpha}{2}|\zeta|^2}d\mu(\zeta)\nonumber\\
&\lesssim \int_{\CC} \int_{D(\zeta,r)} (|w|+1)^q| k_w(z)|^q e^{-\frac{q\alpha}{2}|z|^2}dm(z) d\mu(\zeta)\nonumber\\
&= \int_{\CC}(|w|+1)^q | k_w(z)|^q e^{-\frac{q\alpha}{2}|z|^2}
\int_{D(z, r)}  d\mu(\zeta )dm(z)\nonumber\\
&\leq \int_{\CC} (|w|+1)^q|k_w(z)|^q
e^{-\frac{q\alpha}{2}|z|^2}\mu(D(z,
r))dm(z)= \widetilde{f}_q(w),
\label{spil}
\end{align} where we set  $f(z)=\mu(D(z, r))$. This along with Lemma~\ref{gen} yield
\begin{equation}
\label{normest3}
\|\widetilde{\mu}_q\|_{L^p} \lesssim \|\widetilde{f}_q\|_{L^p}\lesssim \|f\|_{L_{\varphi_q}^p}=\|D_{rq\mu}\|_{L^p}  <\infty
\end{equation} for all $p\geq 1.$  To this end, we remain  with the case for  $p<1.$  Observe that for each $z$ in $ D(z_j,r)$, we  estimate
 \begin{align*}
 |z-w|^2 \geq \big( |w-z_j|-|z-z_j|\big)^2 \geq  |w-z_j|^2-2r |w-z_j|.
 \end{align*}
From this fact  and completing the square in the inner product
from the kernel function, it follows that
\begin{align}
\int_{\CC} |\widetilde{\mu}_q (w)|^p dm(w)&= \int_{\CC} \Bigg(\int_{\CC} (|w|+1)^q |k_w(z)|^q e^{-\frac{\alpha q}{2}|z|^2} d\mu(z)\Bigg)^p d m(w)\ \ \ \ \ \ \  \ \ \ \  \ \ \ \ \ \ \ \ \ \ \ \ \ \ \ \ \ \  \nonumber\\
      &\lesssim \int_{\CC}\Bigg(\sum_{j=1}^\infty \int_{D(z_j,r)} (|w|+1)^q |k_z(w)|^q e^{-\frac{\alpha q}{2}|w|^2}d\mu(z)\Bigg)^p dm(w) \nonumber\\
      &\leq\int_{\CC}\Bigg(\sum_{j=1}^\infty \int_{D(z_j,r)} (|w|+1)^q  e^{-\frac{q\alpha}{2}|w-z_j|^2+r\alpha q |w-z_j|}d\mu(z)\Bigg)^p dm(w)\nonumber\\
        &\le \sum_{j=1}^\infty  \big(\mu(D(z_j,r))\big)^p  \int_{\CC} (|w|+1)^{pq} \frac{ e^{r\alpha qp |w-z_j|}}{e^{-\frac{q p\alpha}{2}|w-z_j|^2}} dm(w),
        \label{otherp}
      \end{align}where  the last inequality is possible since $0<p<1.$
     Applying the change of variables $\zeta=w-z_j$ and using the obvious inequality $|a+b|^q \leq 2^q(|a|^q+ |b|^q)$, the integral above can be estimated as
       \begin{align*}
       \int_{\CC} (|w|+1)^{pq}  e^{-\frac{q p\alpha}{2}|w-z_j|^2+r\alpha qp |w-z_j|} dm(w)\ \ \ \ \ \ \ \ \ \ \ \ \ \ \ \ \  \nonumber\\
       =  \int_{\CC} (|\zeta+z_j|+1)^{pq}  e^{-\frac{q p\alpha}{2}|\zeta|^2+r\alpha qp |\zeta|} dm(\zeta)
       \end{align*} which is bounded by
        \begin{align*}
        2^{qp}(|z_j|+1)^{pq} \int_{\CC} (|\zeta|+1)^{pq}  e^{-\frac{q p\alpha}{2}|\zeta|^2+r\alpha qp |\zeta|} dm(\zeta)
       \lesssim (|z_j|+1)^{pq},
       \end{align*} as the last  integral is easily seen to be convergent. Using this fact, relation \eqref{otherp},
       and the equivalency of the already proved  statements  in (ii) and (iii) in  the lemma, we obtain
        \begin{align}
        \label{normest4}
       \int_{\CC} |\widetilde{\mu}_q (w)|^p dm(w) \lesssim
       \sum_{j=1}^{\infty}(|z_j|+1)^{pq}  \big(\mu(D(z_j,r))\big)^p\simeq \|D_{rq\mu}\|_{L^p}^p
       \end{align} as required. Now from the series of estimates in \eqref{normest1}, \eqref{est1}, \eqref{ess2}, \eqref{normest3} and \eqref{normest4}, we conclude the asymptotic norm relations in \eqref{normre1}.
       \end{proof}

\section{Proof of the main results}
Before we begin proving the main  results, we  proceed to recall a few more facts that will be
used in our subsequent considerations.  A recent result of Constantin \cite{Olivia}   ensures that for
each entire function $f$
 \begin{align}
  \label{olivia}
\int_{\CC} |f(z)|^p e^{-\frac{p\alpha}{2}|z|^2}dm(z) \simeq |f(0)|^p
+\int_{\CC} |f'(z)|^p(1+|z|)^{-p}e^{-\frac{p\alpha}{2}|z|^2}dm(z)
\end{align} for $0<p<\infty$. The corresponding estimate for $p=\infty$ follows from  \cite{TM0} and reads
\begin{align}
\sup_{z\in \CC} |f(z)| e^{-\frac{\alpha}{2}|z|^2}  \simeq |f(0)|+
 \sup_{z\in \CC}  |f'(z)| (1+|z|)^{-1}e^{-\frac{\alpha}{2}|z|^2}.
 \label{tesfa}
\end{align}
These Littelwood--Paley type estimates provide a natural description of  the Fock spaces
$\mathcal{F}_\alpha^p$ in terms of the first derivatives and we
will use them repeatedly  in our proofs. From the relations
in \eqref{olivia} and \eqref{tesfa}, we also easily deduce the
 pointwise estimate
\begin{equation}
\label{point} |f'(z)| \lesssim (1+|z|)e^{\frac{\alpha}{2}|z|^2}
\|f\|_{p}
\end{equation}
for each point $z$ in $\CC$,  exponent $0<p\leq \infty$, and holomorphic function $f.$

\subsection{Proof of Theorem~\ref{thm1} }
(i)
{\it Sufficiency.} First assume that $q<\infty.$
  Since $|f'|^p$ is  subharmonic  for each holomorphic function $f$, by Lemma~1 of \cite{JIKZ}, we have the  local estimate
   \begin{equation}
   \label{subharmonic}
   |f'(z)|^pe^{-\frac{\alpha p}{2}|z|^2} \lesssim \int_{D(z,r)} |f'(w)|^p e^{-\frac{\alpha p}{2}|w|^2}
   dm(w),
   \end{equation} where $r>0$ is fixed. Applying \eqref{olivia} and \eqref{subharmonic} with $p=q$ and
$r=1$, we obtain
\begin{equation}
\label{one} \|J_{(g,\psi)}f\|_{q}^q \lesssim \int_{\CC}
e^{\frac{q\alpha}{2}\big( |\psi(z)|^2-|z|^2\big)}
\frac{|g(z)|^q}{(1+|z|)^q} \int_{\CC}
\chi_{D(\psi(z),1)}(w)\frac{|f'(w)|^q}{e^{\frac{\alpha q}{2}|w|^2}}
dm(w)dm(z).\end{equation} On the other hand, since
$\chi_{D(\psi(z),1)}(w)= \chi_{D(w,1)}(\psi(z))$, by Fubini's
theorem it follows that  the right-hand side of the above
inequality  is equal to
\begin{align}
   \label{II}
   \int_{\CC}|f'(w)|^q&e^{-\frac{\alpha q}{2}|w|^2} \int_{D(w,1)}e^{\frac{q \alpha}{2}|\xi|^2} d\mu_g(\xi)
   dm(w)\nonumber\\
   &\simeq \int_{\CC}|f'(w)|^q\frac{e^{-\frac{\alpha q}{2}|w|^{2}}}{(1+|w|)^q} \int_{D(w,1)}(1+|\xi|)^qe^{\frac{q \alpha}{2}|\xi|^2} d\mu_g(\xi)  dm(w),
      \end{align} where we set $\xi= \psi(z)$,
      \begin{equation*}
d\mu_{g}(E)=\int_{\psi^{-1}(E)} \frac{|g(z)|^q}{(1+|z|)^q}
e^{-\frac{\alpha q}{2}|z|^2} dm(z)
\end{equation*} for every Borel  subset $E$ of $\CC$,
and use the fact that $1+|w| \simeq 1+ |\xi|$ whenever $\xi$
belongs to the disc $D(w,1)$.  Applying \eqref{olivia} with $p=q$
again, the right-hand double integral in \eqref{II} is bounded by
a constant multiple of
\begin{align*}
\|f\|_{q}^q\sup_{w\in \CC} \int_{D(w,1)}(1+|\xi|)^qe^{\frac{q
\alpha}{2}|\xi|^2} d\mu_g(\xi).
\end{align*}
If we show that the supremum above is finite, then the desired
conclusion follows since $\mathcal{F}_\alpha^p \subseteqq
\mathcal{F}_\alpha^q $ whenever $p\leq q.$ To this end, we have
\begin{align}
\sup_{w\in \CC} \int_{D(w,1)}(1+|\xi|)^qe^{\frac{q \alpha}{2}|\xi|^2} d\mu_g(\xi)
&\simeq \sup_{w\in \CC}(1+|w|)^q\int_{D(w,1)}e^{\frac{q \alpha}{2}|\xi|^2} d\mu_g(\xi)\nonumber\\
&\lesssim \sup_{w\in \CC}B_{(|g|^q,\psi)}(w),\label{est0}
\end{align}
where in the last relationship we have used a simple fact that if
$\xi\in D(w,1)$, then
\begin{align}
\label{Ker} |k_w(\xi)|^q=|e^{-\frac{\alpha
}{2}|w|^2+\alpha\overline{w}\xi}|^q=
e^{\frac{q\alpha}{2}(|\xi|^2-|\xi-w|^2)} \gtrsim e^{\frac{q
\alpha}{2}|\xi|^2},
\end{align}
and integrating \eqref{Ker} against the measure $\mu_{g}$ we have
that
\begin{align*}
\int_{D(w,1)} e^{\frac{q \alpha}{2}|\xi|^2}d\mu_g(\xi)  \lesssim
\int_{\CC}|k_w(\xi)|^qd\mu_g(\xi)=\frac{B_{(|g|^q,\psi)}(w)}{(1+|w|)^q}.
\end{align*}

On the other hand, for  $q= \infty,$  the sufficiency of the
condition follows from a simple estimation along with
\eqref{tesfa}  and  \eqref{point} as
 \begin{align}
 \|J_{(g,\psi)} f\|_{\infty}&\simeq \sup_{z\in \CC} \frac{|f'(\psi(z))||g(z)|}{ (1+|z|)}e^{-\frac{\alpha}{2}|z|^2}\nonumber\\
 &\lesssim \|f\|_p  \sup_{z\in \CC} \frac{(|\psi(z)|+1)|g(z)|}{1+ |z|}e^{\frac{\alpha}{2}(|\psi(z)|^2-|z|^2)}
 = \|f\|_p \|M_{(g,\psi)}\|_{L^\infty}.\label{jf}
 \end{align} From the series of   estimations \eqref{one}, \eqref{II}, \eqref{est0} and \eqref{jf}, we have already  one side of the asymptotic relation in
 \eqref{bdd1}, namely that $ \|J_{(g,\psi)} \|$ is bounded by a constant multiple of $\|M_{(g,\psi)}\|_{L^\infty}$ for $q= \infty $ and
 $\|B_{(|g|^q,\psi)}\|_{L^\infty}^{1/q}$ whenever $q<\infty$.

{\it Necessity.} Assuming that $q<\infty,$ we apply $J_{(g,\psi)}$
 to $k_w $  and invoke \eqref{olivia} to obtain
 \begin{equation}
 \label{part}
 1 \gtrsim\|J_{(g,\psi)} k_w\|_{q}^q  \simeq \int_{\CC} \big| wk_w(\psi(z))\big|^q \frac{|g(z)|^q}{(1+|z|)^q}
  e^{-\frac{\alpha q}{2}|z|^2} dm(z)=S_1.
\end{equation} To arrive at the desired conclusion, it remains to show  the estimate
 \begin{align}
 \label{ess4}
 \|J_{(g,\psi)} k_w\|_{q}^q  \gtrsim \int_{\CC} \big|k_w(\psi(z))\big|^q \frac{|g(z)|^q}{(1+|z|)^q}
  e^{-\frac{\alpha q}{2}|z|^2} dm(z)=S_2.
 \end{align}
  To do so, we may  invoke Lemma~\ref{basic}  and equivalently  express the conditions in the theorem
  in terms of the sequence $\big( (1+|z_j|)^q \mu(D(z_j,r))\big)_{j\in\NN}$ for some measure $\mu$ on $\CC.$  Indeed, carefully   examining the proof of the lemma we see that the sequence  $\big( |z_j|^q \mu(D(z_j,r)\big))_{j\in\NN}$
  belongs to $\ell^p$ if and only if the function
  \begin{equation*}
  \widehat{\mu}_q(w)= \int_{\CC} |w|^q |k_w(z)|^q e^{-\frac{\alpha q}{2}|z|^2} d\mu(z)
  \end{equation*} belongs to $L^p.$  Using this fact for the case $p=\infty$, we may equivalently write condition \eqref{part} as
  \begin{equation}
  \label{ess5}
  \sup_{j\in\NN} |z_j|^q \mu_{(g,\psi)}(D(z_j,r)) \simeq \|J_{(g,\psi)} k_w\|_{q}^q \lesssim 1,
  \end{equation} where $\mu_{(g,\psi)}$ is a measure on $\CC$ such that
   \begin{align*}
  d\mu_{(g,\psi)}(z)= \frac{ |g(z)|^q}{(1+|z|)^q}e^{\frac{\alpha q}{2}(|\psi(z)|^2-|z|^2)} dm(z)\circ\psi^{-1}(z).
  \end{align*}
  Having singled out this connection with a discrete setting, we will arrive at \eqref{ess4} if we equivalently show that the sequence
   \begin{align*}
   \sup_{j\in \NN}\mu_{(g,\psi)}(D(z_j,r)) <\infty.
   \end{align*}This immediately follows from \eqref{ess5} if $|z_j|\geq 1$ for all $j$ in $\NN.$ On the other hand, since the sequence $(z_j), j \in \NN$ is an $r/2-$ lattice for $\CC,$ $|z_j| \to \infty $ as $j \to \infty.$ This means that the inequality $|z_j| <1$ can  possibly hold for only a finite number of indices $j,$ and hence
   \begin{align*}
   \sup_{j\in \NN, |z_j|<1}\mu_{(g,\psi)}(D(z_j,r)) \lesssim \sup_{j\in \NN, |z_j|\geq1}\mu_{(g,\psi)}(D(z_j,r))\lesssim\|J_{(g,\psi)} k_w\|_{q}^q
   \end{align*} from  which and \eqref{part} we conclude that
\begin{equation}
\|J_{(g,\psi)} k_w\|_{q}^q \gtrsim S_1+S_2 \gtrsim B_{(|g|^q,\psi)}(w).
\end{equation}
On the other hand, if $q= \infty$, then applying \eqref{tesfa} we
have
\begin{align*}
\|J_{(g,\psi)}\| \gtrsim\|J_{(g,\psi)}k_{w}\|_{\infty} \simeq \sup_{z\in \CC}  \frac{|g(z)| }{1+|z| }|\overline{w} k_{w}(\psi(z))|e^{-\frac{\alpha}{2}|z|^2}\\
 \geq \frac{|g(z)\overline{w}| }{1+|z| }\big|e^{\alpha\overline{w}\psi(z)}\big| e^{-\frac{\alpha}{2}(|w|^2+|z|^2)}
 \end{align*}  for all $z$ and $ w$ in  $\CC.$ In particular,  when we set $w= \psi(z)$, we get
 \begin{equation}
\|J_{(g,\psi)}\|\gtrsim  \|J_{(g,\psi)}k_{w}\|_{\infty} \gtrsim
\frac{|\psi(z)||g(z)|}{1+
|z|}e^{\frac{\alpha}{2}(|\psi(z)|^2-|z|^2)}.
 \label{esm}
 \end{equation}
From \eqref{esm} we have that
\begr
e^{\frac{\alpha}{2}(|\psi(z)|^2-|z|^2)}\lesssim
\frac{|z|+1}{M_\infty(g\psi,|z|)}\label{es1}\endr
whenever $g\psi$ is a nonzero function.
 Estimate
\eqref{es1} along with the fact that the integral mean
$M_\infty(g\psi,|z|)$ is a nondecreasing function of $|z|,$
implies that
\begr \limsup_{|z|\to\infty}(|\psi(z)|-|z|)\le
0.\label{es2}\endr
Otherwise, there would be a sequence
$(w_j)_{j\in\NN}$ such that $|w_j|\uparrow \infty$ as
$j\to\infty$, and
$$\limsup_{j\to\infty}(|\psi(w_j)|-|w_j|)=a>0,$$
 from which and
\eqref{es1} we would have that the sequence
$e^{\frac{\alpha}{2}\big(|\psi(w_j)|^2-|w_j|^2\big)}/(|w_j|+1)$,
$j\in\NN$, is unbounded, and contradicts \eqref{es1}.

From \eqref{es2} and the Cauchy inequality we easily obtain that
$\psi$  has the form $\psi(z)= az+b$ with $|a|\leq 1$ and $b=0$
whenever $|a|=1.$ Using this in \eqref{esm} we have that
\begr
\|J_{(g,\psi)}\|\gtrsim \sup_{z\in\CC}\frac{|az+b||g(z)|}{1+
|z|}e^{\frac{\alpha}{2}(|\psi(z)|^2-|z|^2)},
\label{es3}\endr from
which it easily follows that
\begin{equation}
\label{con}
 \|J_{(g,\psi)}\| \gtrsim  \sup_{z\in\CC}\frac{(|\psi(z)|+1)|g(z)|}{1+ |z|}e^{\frac{\alpha}{2}(|\psi(z)|^2-|z|^2)} = M_{(g,\psi)}(z).
 \end{equation}

(ii)  We now assume that  $B_{(|g|^q,\psi)}$ vanishes at infinity
and proceed to  prove that $ J_{(g,\psi)}$ is compact. We consider
a sequence $(f_n)_{n\in\NN}$ of  functions in
$\mathcal{F}_\alpha^p$ such that $\sup_n \|f_n\|_{p}<\infty$ and
$f_n$ converges to zero uniformly on compact subsets of $\CC$ as $
n\to \infty$. Following the arguments made in the proof of the
sufficiency part, for a positive $R$  and a finite exponent $q\geq
p,$ we have
\begin{align*}
\limsup_{n\to \infty} \|J_{(g,\psi)}f_n\|_{q}^q \lesssim&
\limsup_{n\to \infty}\int_{\CC}\frac{|f_n'(w)|^q}{(1+|w|)^q} e^{-\frac{\alpha q}{2}|w|^2}B_{(|g|^q,\psi)}(w)   dm(w) \ \  \ \ \ \ \ \ \ \ \ \ \nonumber\\
\leq& \limsup_{n\to \infty}\int_{|w|\leq R}\frac{|f_n'(w)|^qe^{-\frac{\alpha  q}{2}|w|^2}}{(1+|w|)^q}B_{(|g|^q,\psi)}(w)  dm(w)\nonumber\\
&+ \limsup_{n\to
\infty}\sup_{|w|>R}B_{(|g|^q,\psi)}(w)\int_{|w|>R}\frac{|f_n'(w)|^qe^{-\frac{\alpha
q}{2}|w|^2}}{(1+|w|)^q} dm(w),\nonumber
\end{align*}
for every $R>0$.

Since $B_{(|g|^q,\psi)}$ vanishes at infinity we have that for
every $\ve>0$ there is an $R_0>0$ such that \begr
\sup_{|w|>\rho}B_{(|g|^q,\psi)}(w)<\ve\label{es3}\endr for every
$\rho\ge R_0.$ We may assume that $R_0=R$.

Applying \eqref{olivia} and  the assumption that $\|f_n\|_{p}$ is
uniformly bounded, we obtain
\begin{align}
\label{do} \limsup_{n\to
\infty}\sup_{|w|>R}B_{(|g|^q,\psi)}(w)\int_{|w|>R}\frac{|f_n'(w)|^qe^{-\frac{\alpha
q}{2}|w|^2}}{(1+|w|)^q}  dm(w) \lesssim
\sup_{|w|>R}B_{(|g|^q,\psi)}(w)<\ve.
\end{align}
On the other hand, we have the next estimate
\begin{align}
\int_{|z|\leq R}\frac{|f_n'(w)|^qB_{(|g|^q,\psi)}(w)}{(1+|w|)^qe^{\frac{\alpha
q}{2}|w|^2}}   dm(w) \lesssim \sup_{|w|\leq
R}|f_n'(w)|^q\lesssim \sup_{|w|\leq 2R}|f_n(w)|^q,\label{es5}
\end{align} where we have used the fact that the integral
\begin{align*}
\int_{|w|\leq R}e^{-\frac{\alpha q}{2}|w|^2} B_{(|g|^q,\psi)}(w)
dm(w)
\end{align*}
is finite due to the boundedness of $B_{(|g|^q,\psi)}(w)$. Since $f_n$ converges to zero uniformly on
compact subsets of $\CC,$ taking the $\limsup$ in \eqref{es5} and combining the result with
 \eqref{do} we get
$$\limsup_{n\to \infty}\|J_{(g,\psi)}f_n\|_{q}^q\lesssim\ve.$$
From this and since $\ve$ is an arbitrary positive number we get
$\lim_{n\to \infty}\|J_{(g,\psi)}f_n\|_q=0$, so by Lemma
\ref{complem} the compactness follows.

We  need to conclude  the same when $q= \infty$. For this, we may
modify a common approach used in dealing with the compactness of
operators acting between spaces of holomorphic function with
target space restricted to be a growth space (see, for example,
\cite{TM0,SS,UEKI2}). To this end, we note that the function
$f_0(z)= z$ belongs to $\mathcal{F}_{\alpha}^p$ for all $p>0.$  It
follows that by the boundedness and \eqref{tesfa},
\begin{equation}
\label{00}
\|J_{(g,\psi)} f_0\|_{\infty} \simeq  \sup_{z\in \CC} |g(z)| (1+|z|)^{-1} e^{-\frac{\alpha}{2}|z|^2}<\infty.
\end{equation}
For each positive $\epsilon,$  the necessity of  the condition
implies that there exists a positive $N_1$  such that
$$ M_{(g,\psi)}(z) <\epsilon$$
for all $z\in\CC$ such that $|\psi(z)| > N_1.$ From this along
with \eqref{tesfa} and \eqref{point}, we get
\begin{equation}
\label{tmo} \frac{|g(z)||f'_n(\psi(z))|}{e^{\frac{\alpha}{2}|z|^2} (1+|z|)}
 \lesssim
\|f_n\|_{p}\frac{|g(z)(|\psi(z)|+1)}{ 1+|z|}
e^{\frac{\alpha}{2}|\psi(z)|^2-\frac{\alpha}{2}|z|^2}\lesssim
M_{(g,\psi)}(z)< \epsilon
\end{equation}
for all $z\in\CC$ such that $|\psi(z)| > N_1$ and all $n.$ On the
other hand if $|\psi(z)|\leq N_1,$ then applying \eqref{00} it is
easily seen that
\begin{equation}
\label{tm1} \frac{|g(z)|}{ 1+|z|}| f_n(\psi(z)) |
e^{-\frac{\alpha}{2}|z|^2} \lesssim \sup_{\{z:|\psi(z)|\leq
N_1\}}|f_n(\psi(z))| \to 0
\end{equation} as $n\to \infty$. Then we combine \eqref{tmo} and \eqref{tm1} to arrive at the desired conclusion.\\
To complete the proof of part (ii) of Theorem~\ref{thm1}, it
remains to verify the necessity of the compactness condition.
Since $k_{z}$ is bounded in $\mathcal{F}_{\alpha}^p$ and uniformly
converges to zero on any compact subsets of $\CC$ as $|z| \to
\infty,$ the compactness of $ J_{(g,\psi)}$ implies
\begin{align*}B_{(|g|^q,\psi)}(z) \lesssim\|J_{(g,\psi)} k_{z}\|_{q}^q\to 0
\end{align*}
as $|z| \to \infty$ whenever $q$ is finite.

Now suppose that $q=\infty$ and  further assume that there exists
a sequence of points $(z_j)_{j\in\NN} \subset \CC$ such that
$|\psi(z_j)| \to \infty$ as $j\to \infty.$ If such a sequence does
not exist, then necessity trivially holds. It follows from the
compactness of $J_{(g,\psi)}$ that
\begin{equation}
\limsup_{j\to \infty} M_{(g,\psi)}(z_j)  \lesssim \limsup_{j\to
\infty}\| J_{(g,\psi)} k_{z_j} \|_{\infty}=0,
\end{equation} and completes the proof of part (ii) of Theorem~\ref{thm1}.\\
The  statements in parts  (iii) and (iv) of the theorem  follow from  simple variants of the proofs of the statements in parts  (i) and (ii) respectively. First observe that $(C_{(g,\psi)}f)'(z) = f'(\psi(z)) g(\psi(z)) \psi'(z)$. It means that we only need
to replace the quantity  $g(z)$ by $ g(\psi(z))\psi'(z)$  in all the above arguments and   proceed as
in the preceding parts. Thus we omit the remaining details.
\subsection{Proof of Theorem~\ref{thm2} }
(i). Since the compactness
obviously implies the boundedness through the normal family
argument, we will prove that boundedness implies the
$L^{p/(p-q)}$ and $L^1$ integrability conditions  and this in turn
implies compactness. Proceeding as in the proof of the  first part of
Theorem~\ref{thm1}, we have
       \begin{align}
       \label{strong1}
       \|J_{(g,\psi)}f\|_{q}^q \lesssim \int_{\CC} \frac{|f'(w)|^q e^{-\frac{q\alpha}{2}|w|^2}} {(1+|w|)^q} \int_{D(w,1)} (1+|\xi|)^q e^{\frac{q\alpha}{2}|\xi|^2} d\mu_g(\xi) dm(w)\nonumber\\
       \lesssim \int_{\CC} \frac{|f'(w)|^q e^{-\frac{q\alpha}{2}|w|^2}} {(1+|w|)^q} B_{(|g|^q,\psi)}(w)dm(w).
       \end{align} Since $p>q,$ applying  H\"{o}lder's inequality, the right hand side  quantity is bounded by
       \begin{align}
       \label{strong2}
      \bigg( \int_{\CC} \frac{|f'(w)|^p e^{-\frac{p\alpha}{2}|w|^2}} {(1+|w|)^p}dm(w)\Bigg)^{\frac{q}{p}} \Bigg(\int_{\CC}|B_{(|g|^q,\psi)}(w)|^\frac{p}{p-q} dm(w))\bigg)^{\frac{p-q}{p}} dm(w)\nonumber\\
      \lesssim \|f\|_{p}^q \Bigg(\int_{\CC}|B_{(|g|^q,\psi)}(w)|^\frac{p}{p-q} dm(w))\Bigg)^{\frac{p-q}{p}}
       \end{align}
whenever $p$ is finite. On the other hand if $p= \infty,$ then by
\eqref{tesfa}, it follows that
       \begin{align}
       \label{twoo}
       \|J_{(g,\psi)}f\|_q^q &\lesssim \int_{\CC} \frac{|f'(w)|^q e^{-\frac{q\alpha}{2}|w|^2}} {(1+|w|)^q} B_{(|g|^q,\psi)}(w) dm(w)\nonumber\\
       &\lesssim \|f\|_{\infty}^q  \int_{\CC} B_{(|g|^q,\psi)}(w) dm(w).
       \end{align} From the estimates in \eqref{strong1},  \eqref{strong2} and \eqref{twoo}, we also  have
       \begin{equation}
       \label{onepart}
              \|J_{(g,\psi)}\|^q \lesssim
              \begin{cases}
       \|B_{(|g|^q,\psi)}\|_{L^{\frac{p}{p-q}}}, & p<\infty\\
       \|B_{(|g|^q,\psi)}\|_{L^1}, & p= \infty.
       \end{cases}
\end{equation}
This establishes one part of the asymptotic relation in
\eqref{bdd3}. The remaining part of the estimate will follow from
our next proof of the integrability condition. For this, we
appeal to the atomic decomposition of functions in Fock spaces,
i.e.,  each function in $\mathcal{F}_\alpha^p$  is generated by an
$\ell^p$  sequence as
   \begin{equation}
   f= \sum_{j=1}^\infty c_j k_{z_j}, \ \ \text{and}\ \ \|f\|_{p} \simeq \|(c_j)\|_{\ell^p}.
   \end{equation} This was proved  in  \cite{SJR} for $p\geq 1$ and in \cite{RW} for $0<p<1$.
We first assume that $0<q<\infty$, and if $(r_j(t))_{j\in \NN}$ is  the Rademacher
sequence of functions on $[0,1]$ chosen  as in \cite{DL}, then
Khinchine's inequality yields
\begin{equation}\Bigg(\sum_{j=1}^\infty
|c_jz_j|^2|k_{z_j}(z)|^2\Bigg)^{\frac{q}{2}} \lesssim \int_0^1
\Bigg|\sum_{j=1}^\infty c_jz_jr_j(t)k_{z_j}(z)\Bigg|^qdt.
\label{Khinchine}
\end{equation}
Note that here if the $r_j(t)$'s are chosen as refereed
above, then the sequence  $(c_jr_j(t))_{j\in \NN}$ belongs to $\ell^p$ with
$\|(c_jr_j(t))\|_{\ell^p} = \|(c_j)\|_{\ell^p}$  for all $t$ and
 \begin{equation}
 \label{Rademacher}
 \sum_{j=1}^\infty c_jr_j(t)k_{z_j}(z) \in \mathcal{F}_{\alpha}^p, \ \text{with}\ \
 \bigg\|\sum_{j=1}^\infty c_jr_j(t)k_{z_j}(z) \bigg\|_{p} \simeq \|(c_j)\|_{\ell^p}.
 \end{equation}
Setting as before
\begin{equation}d\mu_{(g,\psi)}(z)= \frac{ |g(z)|^q}{(1+|z|)^q}e^{\frac{\alpha q}{2}(|\psi(z)|^2-|z|^2)} dm(z)\circ\psi^{-1}(z)
\end{equation}
 and making use of  \eqref{Khinchine},  and subsequently Fubini's
theorem, we obtain
\begin{align}
\label{half}
\int_{\CC} \Bigg( \sum_{j=1}^\infty |c_j z_j|^2
|k_{z_j}(z)|^2\Bigg)^{\frac{q}{2}} d\mu_{(g,\psi)}(z)
\lesssim \int_{\CC}
\Bigg(\int_0^1 \bigg| \sum_{j=1}^\infty c_jz_j
r_j(t)k_{z_j}(z)\bigg|^q
dt\Bigg) d\mu_{(g,\psi)}(z)\nonumber\\
=\int_0^1\Bigg(\int_{\CC}\bigg| \sum_{j=1}^\infty c_jz_j
r_j(t)k_{z_j}(z)\bigg|^q d\mu_{(g,\psi)}(z)\Bigg)dt.
\end{align}
Invoking \eqref{olivia} with $p=q$, the double integral above is
asymptotically equal to
\begin{align}
\int_0^1\Big\|J_{(g,\psi)}\sum_{j=1}^\infty c_j
r_j(t)k_{z_j}\Big\|_{q}^q dt\lesssim\|J_{(g,\psi)}\|^q \|(c_j)\|_{\ell^p}^q,
\label{combine}
\end{align} where the estimate follows by the  boundedness assumption and \eqref{Rademacher}.

Now if $q \geq 2,$ then we obviously have
\begin{align}
\label{alt}
\sum_{j=1}^\infty |c_j|^q |z_j|^q \mu_{(g,\psi)}(D(z_j,r)) \lesssim \int_{\CC}\bigg( \sum_{j=1}^\infty |c_j|^2 |z_j|^2\chi_{D(z_j,r)}(z)\bigg)^{\frac{q}{2}} d\mu_{(g,\psi)}(z).
\end{align} On the other hand, if $q<2$, then   applying H\"older's inequality, we obtain
\begin{align}
\sum_{j=1}^\infty |c_j|^q |z_j|^q \mu_{(g,\psi)}(D(z_j,r)) \leq (N_{\max})^{\frac{2-q}{2}}  \int_{\CC}\Bigg( \sum_{j=1}^\infty |c_j|^2 |z_j|^2\chi_{D(z_j,r)}(z)\Bigg)^{\frac{q}{2}} d\mu_{(g,\psi)}(z)\nonumber\\
 \lesssim  \int_{\CC}\Bigg( \sum_{j=1}^\infty |c_j|^2 |z_j|^2\chi_{D(z_j,r)}(z)\Bigg)^{\frac{q}{2}}
 d\mu_{(g,\psi)}(z).\ \ \ \ \ \ \ \ \label{ses}
\end{align}From \eqref{half}, \eqref{combine}, \eqref{alt},  and \eqref{ses} we deduce
\begin{align}
\label{above}
\sum_{j=1}^\infty |c_j|^q |z_j|^q \mu_{(g,\psi)}(D(z_j,r))\lesssim \|J_{(g,\psi)}\|^q \|(c_j)\|_{\ell^p}^q
\end{align} for each $q.$

Now, if  $p= \infty, $ we set $c_j=1$ for all $j\in \NN$  in
\eqref{above} to see that   the sequence $|z_j|^q \mu(D(z_j,r))
\in \ell^1.$ But we need to show that
\begin{equation}
\label{inf}
\sum_{j=1}^\infty (|z_j|+1)^q \mu_{(g,\psi)}(D(z_j,r)) \lesssim \|J_{(g,\psi)}\|^q.
\end{equation} This  obviously holds when $|z_j|\geq 1$ for all $j\in\NN.$ Thus, we proceed to verify  the
case when $|z_j|<1.$  To this end, note that
since $(z_j)_{j\in \NN}$ is a fixed sequence  with the property that  $|z_j| \to \infty$ as
$j\to \infty$, the inequality $|z_j|<1$ can happen only
for a finite number of $j$'s. It follows that there exist a positive constant $N_f$ for which
\begin{align}
\sum_{j: |z_j|<1}   \frac{\mu_{(g,\psi)}(D(z_j,r))}{(|z_j|+1)^{-q}}\leq
N_f \sum_{j: |z_j|\geq1}  \frac{\mu_{(g,\psi)}(D(z_j,r))}{(|z_j|+1)^{-q}}\lesssim \|J_{(g,\psi)}\|^q.
\end{align}
 On the other hand,  if $ p<\infty, $ then  since $(|c_j|^q)_{j\in\NN} \in \ell^{p/q}$ a duality argument with
 \eqref{above} ensures that  the sequence $  (|z_j|^q \mu(D(z_j,r)))_{j\in\NN}$ belongs to $  \ell^{p/(p-q)}.$
  Using this fact and  following the same arguments made above for the case when $p= \infty,$ we deduce
   \begin{equation}
  \label{fin}
  \big( (|z_j|+1)^q \mu_{(g,\psi)}(D(z_j,r))\big)_{j\in\NN}\in  \ell^{p/(p-q)}.
  \end{equation}
   We combine this with Lemma~\ref{basic} to arrive at the desired conclusion. Looking at the above proof,  we also  have
  \begin{equation*}
  \|J_{(g,\psi)}\|^q \gtrsim \begin{cases}
  \|((|z_j|+1)^q \mu_{(g,\psi)}(D(z_j,r)))\|_{\ell^{p/(p-q)}}, \ \ p<\infty\\
  \|((|z_j|+1)^q \mu_{(g,\psi)}(D(z_j,r)))\|_{\ell^{1}},  \ \ p= \infty
  \end{cases}
  \end{equation*} from which  and \eqref{normre1},  the reverse asymptotic estimate
in \eqref{onepart} holds.

     It remains to prove that $L^{p/(p-q)}$ for $p<\infty$ and $L^1$ for $p=\infty$  integrability  of
      $\widetilde{\mu}_q$        implies compactness of $J_{(g,\psi)}: \mathcal{F}_\alpha^p \to \mathcal{F}_\alpha^q$ whenever $q<p$. To this end, let $f_n$ be a sequence of functions in $\mathcal{F}_\alpha^p$ such  that $\sup_n \|f_n\|_{p}<\infty$ and $f_n$ converges to zero uniformly on compact subset of $\CC$ as $ n\to \infty$. Assume $p<\infty$.  Then  for a positive  $R$, replacing $f$ by $f_n$ in \eqref{strong1} we  write
   \begin{align*}
    \|J_{(g,\psi)}f_n\|_{q}^q \lesssim \Bigg(\int_{|w|\leq R}+ \int_{|w|>R}\Bigg) \frac{|f_n'(w)|^q e^{-\frac{q\alpha}{2}|w|^2}} {(1+|w|)^q} B_{(|g|^q,\psi)}(w)dm(w)\\
    = I_{n1} + I_{n2}. \ \ \ \ \ \
  \end{align*} Applying H\"older's inequality, we estimate the second piece of integral as
  \begin{align*}
  I_{n2}= \int_{|w|>R} \frac{|f_n'(w)|^q e^{-\frac{q\alpha}{2}|w|^2}} {(1+|w|)^q} B_{(|g|^q,\psi)}(w)dm(w)\ \ \ \ \ \ \ \ \ \ \ \ \ \ \ \ \ \ \ \ \
  \ \ \ \ \ \ \  \\
    \leq \|f_n\|_{p}^q \Bigg(\int_{|w|>R}B_{(|g|^q,\psi)}^{\frac{p}{p-q}}(w)dm(w)\Bigg)^{\frac{p-q}{p}}\\
  \lesssim \Bigg(\int_{|w|>R}B_{(|g|^q,\psi)}^{\frac{p}{p-q}}(w)dm(w)\Bigg)^{\frac{p-q}{p}} \to 0
  \end{align*} when $R \to \infty$ since $B_{(|g|^q,\psi)}$ is $L^{p/(p-q)}$ integrable.  On the other hand, because of this integrability and H\"older's inequality again, we have
  \begin{align*}
  I_{n1} = \int_{|w|\leq R} \frac{|f_n'(w)|^q e^{-\frac{q\alpha}{2}|w|^2}} {(1+|w|)^q} B_{(|g|^q,\psi)}(w)dm(w)\lesssim \bigg(\int_{|w|\leq R} \frac{|f_n'(w)|^p e^{-\frac{q\alpha}{2}|w|^2}} {(1+|w|)^p} dm(w)\bigg)^{\frac{q}{p}}.
  \end{align*} Now for sufficiently large $R,$ by \eqref{olivia}, we have that
  \begin{align*}
  I_{n1} \lesssim \Bigg(\int_{|w|\leq R} \frac{|f_n'(w)|^p e^{-\frac{q\alpha}{2}|w|^2}} {(1+|w|)^p} dm(w)\Bigg)^{\frac{q}{p}} \simeq \Bigg(\int_{|w|\leq R} |f_n(w)|^p e^{-\frac{p\alpha}{2}|w|^2} dm(w)\Bigg)^{\frac{q}{p}}\\
   \lesssim \sup_{|w|\leq R} |f_n(w)|^q \to 0
  \end{align*} as $n \to \infty $ since $f_n$ converges to zero uniformly on compact subsets of $\CC.$\\
  Similarity, when $p= \infty$, by  \eqref{point} and  since $\|f_n\|_{\infty}$ is uniformly bounded, it follows that
\begin{align*}
  I_{n2} \leq \|f_n\|_{\infty}^q \int_{|w|>R}B_{(|g|^q,\psi)}^{\frac{p}{p-q}}(w)dm(w)\to 0
  \end{align*} when $R \to \infty$. For such sufficiently big $R,$ we estimate the remaining piece of integral as
  \begin{align*}
  I_{n1} \lesssim \sup_{|w|\leq R}\frac{|f_n'(w)|^q } {(1+|w|)^q} \int_{|w|\leq R} e^{-\frac{q\alpha}{2}|w|^2}B_{(|g|^q,\psi)}(w)dm(w) \lesssim  \sup_{|w|\leq R}|f_n(w)|^q \to 0
        \end{align*} as $n \to \infty $.\\
                The proof of part (ii) of the theorem is very similar to the proof made  above for part (i) and we omit it.
\subsection{Proof of Corollary~\ref{cor1}} (i) As pointed out earlier, for
the special case when  $\psi(z)= z$, the operators $J_{(g,\psi)}$
reduce to the Volterra companion operator $J_g.$   We may first
assume that $g$ is a constant function and $q<\infty$. Then
setting $\psi(z)=z$ in Theorem~\ref{thm1}, we see
   \begin{align*}
   \sup_{w\in\CC}B_{(|g|^q,\psi)}(w)\simeq  \sup_{w\in\CC}\int_{\CC} (|w|+1)^q| k_w(z)|^q \frac{e^{-\frac{q\alpha}{2}|z|^2}}{(1+|z|)^q}dm(z)\\
  = \sup_{w\in\CC} \int_{\CC} (|w|+1)^q   \frac{e^{-\frac{q\alpha}{2}|w-z|^2}}{(1+|z|)^q}dm(z)
      \end{align*} is finite. In a similar way, if $q= \infty,$ then
          $     M_{(g,\psi)}(z)\simeq 1$ for all points $z$
in $\CC.$      On the other hand, if $B_{(|g|^q,\psi)}$ is bounded, then by subharmonicity, we have
   \begin{align}
   \sup_{w\in\CC} B_{(|g|^q,\psi)}(w)\geq \int_{D(w,1)} (|w|+1)^q|k_w(z)|^q \frac{e^{-\frac{q\alpha}{2}|z|^2}|g(z)|^q}{(1+|z|)^q}dm(z)\nonumber\\
   \gtrsim \frac{|wg(w)|^q}{(1+|w|)^q},
   \label{nec}
   \end{align} where we  used the fact that $1+|z| \simeq 1+|w|$ whenever $z$ belongs to the disc $D(w,1)$. The above estimate implies that
\begin{equation*}
\sup_{w\in\CC}|g(w)| <\infty.
\end{equation*} Since $g$ is an entire function, this holds only when  $g$ is a constant function. For $q= \infty$, the necessity is rather immediate because
\begin{equation*}
\infty > \sup_{z\in \CC}M_{(g,\psi)}(z)\gtrsim \sup_{z\in
\CC}|g(z)|.
\end{equation*} We next show the claim that $M_g$ is bounded if and only if $g$ is a constant function again.  Assuming that $M_g$
is bounded and $q<\infty,$ we have
\begin{eqnarray}
\|M_g k_w\|_{q}^q &=& \frac{\alpha q}{2\pi}\int_{\CC} |g(z)|^q |k_w(z)|^q e^{-\frac{q\alpha}{2}|z|^2}dm(z)\nonumber\\
&\gtrsim& \int_{D(w,1)}|g(z)|^q e^{-\frac{q\alpha}{2}|z-w|^2}dm(z)\geq |g(w)|^q
\label{necc}
\end{eqnarray} for all $w$, where the last inequality follows by subharmonicity again. From this we deduce that $g$ is a constant function.\\
Conversely,  assume that $g$ is a constant function. Then  for each $f$ in $\mathcal{F}_\alpha^p,$  we have
$\|M_g f\|_{q}  \simeq \|f\|_{q} \leq \|f\|_{p} $ where we use the inclusion property $\mathcal{F}_\alpha^p \subset \mathcal{F}_\alpha^q $ whenever $p\leq q.$
For the case when  $q= \infty,$ observe that a constant $g$ implies $\|M_g f\|_{\infty}  \simeq \|f\|_{\infty} \leq \|f\|_{p} $ for
any $p$ from which  the boundedness of  $M_g$ follows.  On the other hand, if $M_g$ is bounded, then
\begin{equation*}
\|M_g k_w\|_{\infty}\geq |g(z) k_w(z)|e^{-\frac{\alpha}{2}|z|^2}
\end{equation*} for each $w$ and $z$ in $\CC$. In particular when we  set $w=z,$ we find that $|g(w)|$ is uniformly bounded
independent of $w$ from which the assertion follows.

To prove the corresponding statements for compactness, we may note that if $g=0,$ then both $J_g$ and $M_g$ are the trivial zero maps and they are compact. On the other hand,  suppose  $J_g$ is compact. Then by Theorem~\ref{bdd1} and \eqref{nec}, we have that
 \begin{equation}
 0= \lim_{ |w| \to \infty} B_{(|g|^q,\psi)}(w)\geq \lim_{ |w| \to \infty} |g(w)|.
 \end{equation} This holds only if $g= 0$. The case for $q= \infty$ is straightforward.

  Since $k_w$ is a unit norm functions which converges uniformly to zero on compact subset of $\CC$ as $|w| \to \infty$, by \eqref{necc} we have
\begin{equation*}
0= \lim_{ |w| \to \infty}\|M_g k_w\|_{q} \geq \lim_{ |w| \to \infty} |g(w)|
\end{equation*} from which the desired conclusion follows again.

(ii) Since the sufficiency is trivial, we shall assume that  $J_g$ is bounded (compact) and proceed to show that $g$ is the zero function. We may first suppose $p<\infty.$
Then an application of  part (i) of Theorem~\ref{thm2} and subharmonicity  give
\begin{eqnarray*}
\int_{\CC} B_{(|g|^q,\psi)}^{\frac{p}{p-q}} (w)dm(w)&
\geq& \int_{\CC}\Bigg(\int_{D(w,1)} \frac{(|w|+1)^q |k_w(z)g(z)|^q}{(1+|z|)^qe^{\frac{\alpha q}{2}|z|^2}}dm(z)\Bigg)^{\frac{p}{p-q}} dm(w)\nonumber\\
&\geq&
\int_{\CC} |g(w)|^{\frac{qp}{p-q}} dm(w).
\end{eqnarray*}This holds only if $g=0$. On the other hand, if $p= \infty, $ then we repeat the above algorithm with exponent $p/(p-q)$ replaced by
$1$ and to easily arrive at  the same conclusion.

To prove that bounded (compact) $M_g$ implies $g$ is the zero
function, we act as in the proof of the necessity of the condition
in Theorem~\ref{thm2}. We in particular follow the rout leading to
the estimate \eqref{combine}. In this case, the corresponding
estimate would be
\begin{equation}
\int_{\CC} \Bigg( \sum_{j=1}^\infty |c_j |^2
|k_{z_j}(z)|^2\Bigg)^{\frac{q}{2}} d\mu(z)\lesssim \|M_{g}\|^q
\|(c_j)\|_{\ell^p}^q, \label{combined}
\end{equation} where $d\mu_g(z)=|g(z)|^qdm(z)$.
On the other hand,
\begin{align}
\int_{\CC} \Bigg( \sum_{j=1}^\infty |c_j |^2
|k_{z_j}(z)|^2\Bigg)^{\frac{q}{2}} d\mu_g(z) \geq \frac{1}{N_{\max}} \sum_{j=1}^\infty \int_{D(z_j,2r)}|c_j k_{z_j}(z)g(z)|^q e^{-\frac{\alpha q}{2}|z|^2} dm(z)\nonumber\\
\gtrsim \sum_{j=1}^\infty \int_{D(z_j,2r)}|c_j |^q |g(z)|^q  dm(z), \nonumber
\end{align} where the first inequality follows by Lemma~\ref{covering}.  Since the sequence $(c_j)_{j\in \NN}$ was chosen arbitrarily from $\ell^p$, the above relation together with \eqref{combined} implies the sequence
\begin{equation*}
\int_{D(z_j,2r)} |g(z)|^q  dm(z)
\end{equation*} belongs to $ l^{p/(p-q)}$ for $p<\infty$ and $l^1$  for $p=\infty$. By the subharmonicity, we have
\begin{align}
\label{subharm}
|g(w)|^q \lesssim \int_{D(z_j,2r)} |g(z)|^q  dm(z)
\end{align}  for each $ w$ in the disc $D(z_j,3r/2)$ and  for each $i\geq 1$.  From this it follows
for finite $p$ that
\begin{align}
\label{finite}
\int_{\CC} |g(w)|^{\frac{pq}{p-q}} dm(w)\lesssim \sum_{j=1}^\infty \Bigg(\int_{D(z_j,3r/2)} |g(z)|^q  dm(z)\Bigg)^{\frac{p}{p-q}} \ \ \ \ \ \ \ \ \  \ \ \ \ \ \ \ \ \ \nonumber\\
\lesssim\sum_{j=1}^\infty \Bigg(\int_{D(z_j,r)} |g(z)|^q  dm(z)\Bigg)^{\frac{p}{p-q}}<\infty.
\end{align} Seemingly, if $p= \infty,$ then \eqref{subharm} ensures
\begin{equation}
\label{infinite}
\int_{\CC} |g(w)|^{q}dm(w) \lesssim \sum_{j=1}^\infty \int_{D(z_j,r)} |g(z)|^q  dm(z)<\infty.
\end{equation} Since  $g$ is analytic,
the estimates in \eqref{finite}  and \eqref{infinite} hold only if
 $g$ is the zero function as asserted. Interested readers may consult \cite{RA}  to see  why the zero function
 is the only $L^q$ integrable  entire function on $\CC$.\\


\begin{thebibliography}{BRSHZE}
\bibitem{Alman} A. Aleman, A class of integral operators on spaces of analytic functions, {\it Topics in
complex analysis and aperator theory,} 3--30, Univ. M\'{a}laga,
M\'{a}laga, 2007.

\bibitem{ALC} A. Aleman and J.  Cima, An integral operator on $H^p$ and Hardy's inequality, {\it J. Anal.
Math.} \textbf{85}  (2001), 157--176.

\bibitem{Alsi1} A. Aleman  and A. Siskakis, An integral operator on $H^p$, {\it Complex Variables} \textbf{28} (1995), 149--158.

\bibitem{Alsi2} A. Aleman  and A. Siskakis, Integration operators on Bergman spaces, {\it Indiana University Math J.} \textbf{46} (1997), 337--356.

\bibitem{Olivia} O. Constantin, Volterra type integration operators on Fock spaces, {\it Proc. Amer. Math. Soc.} \textbf{140} (12) (2012), 4247--4257.

\bibitem{FJ} R.  Fleming and J. Jamison, \emph{Isometries on Banach spaces, function spaces,
monographs and surveys in pure and applied mathematics}, vol. 129. Chapman and Hall/CRC, Boca Raton (2003).

\bibitem{ZHXL} Z. Hu and X. Lv, Toeplitz operators from one Fock space to
another, {\it Integr. Equ. Oper. Theory} \textbf{70} (2011),
541--559.
\bibitem{JIKZ} J. Isralowitz and  K. Zhu, Toeplitz operators on the Fock
space, \emph{Integr. Equ. Oper. Theory},  \textbf{66} (2010), no. 4, 593--611.

\bibitem{SJR} S. Janson, J. Peetre, and R. Rochberg, Hankel forms and the Fock space,
{\it Rev. Mat. Iberoamericana} \textbf{3} (1987), 61--138.

\bibitem{SLI} S.~Li, Volterra composition operators between
weighted Bergman spaces and Bloch type spaces, {\it J. Korean
Math. Soc.} \textbf{45} (2008), 229--248.


\bibitem{LIS} S.~Li and S.~Stevi\' c, Generalized composition
operators on Zygmund spaces and Bloch type spaces, {\it J. Math.
Anal. Appl.} {\bf 338} (2008), 1282--1295.

\bibitem{jmaa345} S.~Li and S.~Stevi\' c, Products of Volterra type
operator and composition operator from $H^\infty$ and Bloch spaces
to the Zygmund space, {\it J. Math. Anal. Appl.} {\bf 345} (2008),
40--52.

\bibitem{DL} D.~Luecking, Embedding theorems for space of analytic functions via Khinchine's inequality, {\it Michigan Math. J.}
\textbf{40} (1993),  333--358.

\bibitem{TM} T. Mengestie, \textit{Product of Volterra type integral and composition operators on weighted Fock spaces},
 J. Geom. Anal., \textbf{24}(2014), 740--755.

\bibitem{TM0} T. Mengestie, \textit{Volterra type and weighted composition operators on  weighted Fock spaces},  Integr. Equ. Oper.
Theory, \textbf{76} (2013), no 1, 81--94.

\bibitem{Pom} C. Pommerenke, Schlichte Funktionen und analytische Funktionen von beschr\'{a}nkter mittlerer
Oszillation, {\it Commentarii Mathematici Helvetici} \textbf{52}
(1977), no. 4, 591--602.

 \bibitem {RA} Rashkovskii, A. Rashkovskii, Classical and new loglog-theorems, {\it Expo. Math.,} \textbf{27} (2009), no. 4, 271--287.

     \bibitem{KH} K. Seip and  El. Youssfi, Hankel operators on Fock spaces and related Bergman kernel estimates,  {\it Journal of Geometric Anal}., \textbf{23}(2013), 170--201.

\bibitem{Si} A. Siskakis, Volterra operators on spaces of analytic functions-a survey, {\it Proceedings of the first advanced course
in operator theory and complex analysis, 51--68, Univ. Sevilla
Secr., Seville,} 2006.


\bibitem{SS} S. Stevi\'c, Weighted composition operators between Fock-type spaces in
${\mathbb C}^N$, {\it Appl. Math. Comput.} {\bf 215} (2009),
2750--2760.



\bibitem{UEKI2} S. I. Ueki, Weighted composition operator on some
function spaces of entire functions, {\it Bull. Belg. Math. So.
Simon Stevin} \textbf{17} (2010), 343--353.

\bibitem{UK2} S. I. Ueki, On the Li--Stevi\'c integral type operators from weighted Bergman spaces into $\a$-Zygmund spaces,
{\it Integral Equations Operator Theory}  \textbf{74} (1) (2012),
137--150.

\bibitem{RW} R. Wallst$\acute{e}$n, The $S^p$ Criterion for Hankel forms on the Fock space, $0<p<1$, {\it Math. Scand.} \textbf{64} (1989), 123--132.
\bibitem{KZH} K. Zhu, Spaces of Holomorphic Functions in the Unit Ball, Springer-Verlag, New
York, 2005.

\end{thebibliography}
\end{document}